\documentclass[12pt]{article}

\usepackage{amstext    }
\usepackage{amsthm    }
\usepackage{a4}
\usepackage[mathscr]{eucal}
\usepackage{mathrsfs}

\usepackage{showkeys}

\usepackage{amsmath}
\usepackage{amssymb}
\usepackage{amscd}

\newtheorem{theorem}{Theorem}[section]

\newtheorem{corollary}[theorem]{Corollary}
\newtheorem{lemma}[theorem]{Lemma}
\newtheorem{remark}[theorem]{Remark}

\newcommand{\SO}{{\rm SO}}
\newcommand{\PSL}{{\rm PSL}}

\newcommand{\R}{\mathbb{R}}

\title{Equidistribution of periodic points for modular correspondences}

\author{Tien-Cuong Dinh}

\begin{document}

\maketitle

\begin{abstract}
Let $T$ be an exterior modular correspondence on an irreducible locally symmetric space $X$.
In this note, we show that the isolated fixed points of the power $T^n$ are equidistributed with respect to the invariant measure on $X$ as $n$ tends to infinity. A similar statement is given for general sequences of modular correspondences.
\end{abstract}

\noindent
{\bf Classification AMS 2010:} 37A45, 37A05, 11F32


\section{Introduction} \label{intro} 

Let $G$ be a connected Lie group and  $\Gamma\subset G$ be 
a torsion-free lattice.
Let $\widehat\lambda$ denote the probability measure on $\widehat X:=\Gamma\backslash G$ induced by the invariant measure on $G$. Consider also an element $g\in G$ such that $g^{-1}\Gamma g$ is commensurable with $\Gamma$, that is, $\Gamma_g:=g^{-1}\Gamma g\cap \Gamma$ has finite index in $\Gamma$. Denote by $d_g$ this index.

The map $x\mapsto (x,gx)$ induces a map from $\Gamma_g\backslash G$ to $\widehat X\times\widehat X$. Let $\widehat Y_g$ be its image. The natural projections $\widehat\pi_1,\widehat\pi_2$ from $\widehat Y_g$ onto the factors of $\widehat X\times\widehat X$ define two coverings of degree $d_g$. Both of them are Riemannian with respect to every left-invariant Riemannian metric 
on $G$. The correspondence $\widehat T_g$ on  $\widehat X$ associated with $\widehat Y_g$ is called {\it irreducible modular}. 

{\it A general modular correspondence} $\widehat T$ on $\widehat X$ is a finite sum of irreducible ones, i.e. $\widehat T$ is associated with a sum $\widehat Y=\widehat Y_{g_1}+\cdots+\widehat Y_{g_m}$ that we call {\it the graph} of $\widehat T$. The degree $d$ of $\widehat T$ is the sum of the degrees of $\widehat T_{g_i}$. We refer the reader to \cite{ClozelOtal, ClozelUllmo, Margulis} for more details.

If $a$ is a point in $\widehat X$, define 
$\widehat T(a):=\widehat \pi_2(\widehat\pi_1^{-1}(a))$ and $\widehat T^{-1}(a):=\widehat\pi_1(\widehat\pi_2^{-1}(a))$. They are sums of $d$ points which are not necessarily distinct. If $\widehat U$ is a small neighbourhood of $a$, the restriction of $\widehat T$ to $\widehat U$ can be identified to $d$ local isometries $\widehat \tau_i:\widehat U\to \widehat U_i$  from $\widehat U$ to neighbourhoods $\widehat U_i$ of points $a_i$ in $\widehat T(a)$. 
All these isometries are induced by left-multiplication by elements of $G$. 
If $a$ is a fixed point of $\widehat\tau_i$, i.e. $a=a_i$, we say that $a$ is {\it a fixed point}  
of $\widehat T$. 
When $a$ is an isolated fixed point of $\tau_i$ we also say 
$a$ is {\it an isolated fixed point} of $T$. These points are repeated according to their multiplicities.

The composition $\widehat T\circ \widehat S$ of two modular correspondences $\widehat T$ and $\widehat S$ can be obtained by composing the above local isometries. This is also a modular correspondence. Its degree is equal to $\deg(\widehat T)\deg(\widehat S)$. Even when $\widehat T$ and $\widehat S$ are irreducible, their composition is not always irreducible. Denote by $\widehat T^n:=\widehat T\circ \cdots \circ \widehat T$, $n$ times, {\it the iterate of order $n$} of $\widehat T$. {\it Periodic points of order $n$} of $\widehat T$ are fixed points of $\widehat T^n$. 

Let $\mu$ be a probability measure on $\widehat X$. Define a positive measure $\widehat T_*(\mu)$ of mass $d$ on $\widehat X$ by
$$\widehat T_*(\mu):=(\widehat\pi_2)_* (\widehat\pi_1)^*(\mu).$$
A sequence of correspondences $\widehat T_n$ of degree $d_n$ is said to be {\it equidistributed} if for any $a\in\widehat X$ the sequence of probability measures $d_n^{-1} (\widehat T_n)_*(\delta_a)$ converges weakly to $\widehat\lambda$ as $n$ tends to infinity. Here, $\delta_a$ denotes the Dirac mass at $a$.

Let $K$ be a compact Lie subgroup of $G$. 
Since the left-multiplication on $G$ commutes with the right-multiplication, a modular correspondence $\widehat T$ as above, induces a modular correspondence $T$ on $X:=\widehat X/K$ with the same degree. Its graph is the projection $Y$ of $\widehat Y$ on $X\times X$. The above notion and description of $\widehat T$ can be extended to $T$ without difficulty. We call $\widehat T$  {\it the lift} of $T$ to $\widehat X$. 
Consider on $X$ the probability measure $\lambda$ induced by the invariant measure on $G$, i.e. the direct image of $\widehat\lambda$ in $X$.
Here is our main result. 

\begin{theorem} \label{th_main}
Let $T_n$ be a sequence of modular correspondences on $X$ and let $\widehat T_n$ be the lifts of $T_n$ to $\widehat X$. Assume that the sequence $\widehat T_n$ is equidistributed. Then the isolated fixed points of $T_n$ are equidistributed. More precisely, there is a constant $s\geq 0$, depending only on $G$ and $K$, such that if $d_n$ is the degree of $T_n$ and $P_n$ is the set of isolated fixed points of $T_n$ counted with multiplicity, we have
$$\lim_{n\to\infty} {1\over d_n}\sum_{a\in P_n}\delta_a=s\lambda.$$
\end{theorem}

The last convergence is equivalent to the following property. If $W$ is an open subset of $X$ such that its boundary has zero $\lambda$ measure, then
$$\lim_{n\to\infty} {|P_n\cap W| \over d_n}=s\lambda(W).$$
We can of course replace $W$ with $\overline W$. 

Now, assume moreover that $G$ is semi-simple, $K$ is a maximal compact Lie subgroup of $G$ and $\Gamma$ is an irreducible lattice.  An irreducible correspondence $T$ associated with an element $g\in G$ as above is {\it exterior} 
if the group generated by $g$ and $\Gamma$ is dense in $G$. For such a correspondence, Clozel-Otal proved in \cite{ClozelOtal} that the iterate sequence $\widehat T^n$ is equidistributed (their proof given for $T$ is also valid for $\widehat T$), see also Clozel-Ullmo \cite{ClozelUllmo}. We deduce from Theorem \ref{th_main} the following result.

\begin{corollary}
Let $T$ be an exterior correspondence on an irreducible locally symmetric space $X$ as above. Then the isolated periodic points of order $n$ of $T$ are equidistributed with respect to $\lambda$ as $n$ tends to infinity. 
\end{corollary}

The proof of our main result will be given in Section \ref{section_proof}. In Section \ref{section_rk}, we will give 
similar results related to the Arnold-Krylov-Guivarc'h theorem \cite{ArnoldKrylov, Guivarch}.
We refer to Benoist-Oh \cite{BenoistOh} and Clozel-Oh-Ullmo \cite{ClozelOhUllmo} for other sequences of modular correspondences for which our main result can be applied. The reader will also find in  Clozel-Ullmo \cite{ClozelUllmo}, Dinh-Sibony \cite{Dinh, DinhSibony}  and Mok-Ng \cite{Mok, MokNg} some related topics. 

\bigskip
\noindent
{\bf Acknowledgment.} The author would like to thank Professor Nessim Sibony for help during the preparation of this paper. 

\section{Proof of the main result} \label{section_proof}

Fix a Riemannian metric on $G$ which is invariant under the left-action of $G$ and the right-action of $K$. 
It induces Riemannian metrics on $\widehat X$ and $X$. 
We normalize the metric so that the associated volume form on $X$ is a probability measure. So, it is equal to 
$\lambda$. If $\Pi:\widehat X\to X$ is the canonical projection, we have $\Pi_*(\widehat\lambda)=\lambda$. Let $l$ and $m$ denote the dimension of $G$ and $X$ respectively. 

Fix a point $c\in X$. Denote by $B(c,r)$ the ball of center $c$ and of radius $r$ in $X$. In order to prove the main result, we will consider
the following quantity 
$${|P_n\cap B(c,r)| \over d_n}\cdot$$

Let $\Phi$ denote the natural projection from $G$ to $\widetilde X:=G/K$. 
The image of $K$ by $\Phi$ is a point that we denote by 0. Denote by $B(0,r)$ the ball of center 0 and of radius $r$ in $\widetilde X$. Define $K_r:=\Phi^{-1}(B(0,r))$.  
So, $K_r$ is a union of classes $xK$ with $x\in G$.
Fix also a constant $r_0>0$ small enough so that $B(0,r')$ is convex for every $r'\leq 3r_0$.  
Here, the convexity is with respect to the Riemannian metric induced by the one on $G$. 
From now on, assume that $r<r_0$.

\begin{lemma} \label{lemma_char_fix}
Let $g$ be an element of $G$.
If $g$ admits a fixed point in $B(0,r)$ then $g$ belongs to $K_{2r}$. 
The set of fixed points of $g$ in $B(0,r_0)$ 
is a convex submanifold of $B(0,r_0)$. Moreover, a fixed point $e\in B(0,r_0)$ of $g$ is isolated if and only if
$1$ is not an eigenvalue of the differential of $g$ at $e$.
\end{lemma}
\proof
Assume that $g$ admits a fixed point $e$ in $B(0,r)$.
Since $g$ is locally isometric, $g(0)$ belongs to $B(0,2r)$. It follows that $g$ belongs to $K_{2r}$. 
If $e,e'$ are two different fixed points in $B(0,r_0)$ then every point of the geodesic in $B(0,r_0)$ containing $e,e'$ is fixed. We deduce that the set of fixed points in $B(0,r_0)$ is a convex submanifold. If 1 is an eigenvalue of the differential of $g$ at $e$,  the associated tangent vector at $e$ defines a geodesic of fixed points. This implies the last assertion in the lemma.
\endproof

Recall that {\it a semi-analytic set} in a real analytic manifold $W$ is locally defined by a finite family of inequalities $f>0$ or $f\geq 0$ with $f$ real analytic. A set in $W$ is {\it subanalytic} if locally it is the projection on $W$
of a bounded semi-analytic set in $W\times\R^n$. The boundary of a subanalytic open set is also subanalytic with smaller dimension. We refer the reader to \cite{Lojasiewicz} for further details.
We will need the following lemma.

\begin{lemma} \label{lemma_subanalytic}
Let $M_r$ denote the set of all $g\in G$ which admit exactly one fixed point in $B(0,r)$. 
Then $M_r$ is a subanalytic open set contained in $K_{2r}$.
\end{lemma}
\proof 
The last assertion in Lemma \ref{lemma_char_fix} implies that $M_r$ is open. The first assertion of this lemma implies that $M_r$ is contained in $K_{2r}$. 

Denote by $M'$ the set of points $(g,x)$ in $K_{2r_0}\times B(0,r_0)$ such that $g(x)=x$. This is an analytic subset of $K_{2r_0}\times B(0,r_0)$. So, it is a semi-analytic set in $G\times \widetilde X$. 
Let  $M$ be the set of points $(g,x)$ in $M'$ such that the differential of $g$ at $x$ does not have 1 as eigenvalue. So, $M$ is also a semi-analytic set. 
 
If $\sigma_1,\sigma_2$ are the natural projections from $M'$ to $G$ and to $\widetilde X$ respectively,  we deduce from Lemma \ref{lemma_char_fix} that $M_r$ is equal to $\sigma_1(M\cap\sigma_2^{-1}(B(0,r)))$. Moreover, $\sigma_1$ defines a bijection from $M\cap\sigma_2^{-1}(B(0,r))$ 
to $M_r$. It is now clear that $M_r$ is a subanalytic set. 
\endproof

Consider a general modular correspondence $T$ as above. Let $\pi_1,\pi_2$ denote the natural projections from $Y$ to $X$. If $r$ is small enough, the ball $B(c,r)$ is simply connected and 
$\pi_1^{-1}(B(c,r))$ is the union of $d$ balls 
$B(c'_i,r)$ of center $c'_i$ in $Y$. The restriction of $\pi_1$ to $B(c'_i,r)$ is injective. The projection $\pi_2$ sends $B(c'_i,r)$ to the ball $B(c_i,r)$ of center $c_i:=\pi_2(c_i')$ in $X$.  So, the restriction of $T$ to $B(c,r)$ is identified with the family of $d$ maps $\tau_i:B(c,r)\to B(c_i,r)$.

Fix a point $b\in \widehat X$ such that $\Pi(b)=c$. Let $\widehat T$ denote the lift of $T$ to $\widehat X$ as above. The restriction of $\widehat T$ to $B(b,r)$ can be identified with a family of $d$ maps $\widehat\tau_i:B(b,r)\to B(b_i,r)$ which are the lifts of $\tau_i$ to $\widehat X$, i.e. we have $\Pi\circ\widehat\tau_i=\tau_i\circ \Pi$. 
 
Fix also a point $a\in G$ such that $\Psi(a)=b$ where $\Psi:G\to \widehat X$ is the natural projection. The left-multiplication by $a$ induces the map $x\mapsto \Psi(ax)$ from $M_r$ to $\widehat X$. Its image is independent of the choice of $a$ and is denoted by $M_{b,r}$. Since $\Gamma$ is torsion-free, its intersection with $K$ is trivial. Therefore, when $r$ is small enough, the above map is injective on $K_{2r}$. So, it defines a bijection from $M_r$ to $M_{b,r}$. This is an isometry since the metric on $G$ is invariant.

\begin{lemma} \label{lemma_fix_bis}
The map $\tau_i$ admits exactly one fixed point in $B(c,r)$ if and only if $b_i$ belongs to $M_{b,r}$. 
\end{lemma}
\proof
Without loss of generality, we can assume that $T$ and $\widehat T$ are irreducible and given by an element $g\in G$ such that $g^{-1}\Gamma g$ is commensurable with $\Gamma$. Choose $d$ elements $\delta_1,\ldots,\delta_d$ of $\Gamma$ which represent the classes of $\Gamma_g\backslash\Gamma$. Then, up to a permutation,  $\widehat\tau_i$ and $\tau_i$ are induced by the maps $x\mapsto g_ix$ where $g_i:=g\delta_i$. 

Assume that $\tau_i$ has a unique fixed point in $B(c,r)$. This point can be written as $\Theta(ae)$ for some 
point $e\in B(0,r)$, where $\Theta$ is the canonical projection from $\widetilde X$ to $X$. So, we have $g_iae=\gamma ae$ for some $\gamma\in \Gamma$. The maps $\widehat\tau_i$ and $\widehat \tau$ are also induced by $x\mapsto g_i'x$ where $g_i':=\gamma^{-1}g_i$ since $\gamma^{-1}\in\Gamma$. We have $g_i'ae=ae$ and $(a^{-1}g_i'a)e=e$. By Lemma \ref{lemma_char_fix}, $a^{-1}g_i'a$ belongs to $M_r$.  Since $b_i=\Psi(g_i'a)$, we deduce that
$b_i\in \Psi(aM_r)=M_{b,r}$. We see in the above arguments that the converse is also true.
\endproof

\noindent
{\bf End of the proof of Theorem \ref{th_main}.}
Denote by $\lambda'$ the volume form on $G$ which induces on $\widehat X$ the form $\widehat\lambda$.
By Lemma \ref{lemma_subanalytic}, $M_r$ and $M_{b,r}$ are subanalytic sets. So, their boundaries are of dimension $\leq l-1$. 
Since the sequence $\widehat T_n$ is equidistributed, using Lemma \ref{lemma_fix_bis}, we obtain
$$\lim_{n\to\infty}{|P_n\cap B(c,r)| \over d_n}=\lim_{n\to\infty} {|\widehat T_n(b)\cap M_{b,r}| \over d_n}=\widehat\lambda(M_{b,r})=\lambda'(M_r).$$

It follows that the sequence of positive measures
$${1\over d_n} \sum_{x\in P_n} \delta_x$$
converges to a measure $\mu$ which satisfies $\mu(B(c,r))=\lambda'(M_r)$ for $r$ small enough.
Since $M_r$ is contained in $K_{2r}$, the last quantity is of order $O(r^m)$. Hence, $\mu= s\lambda$ where $s\geq 0$ is a function. Finally, the fact that $\lambda'(M_r)$ is independent of $c$ implies that $s$ is constant. It depends only on $G$ and $K$.
\hfill $\square$

\begin{remark} \rm
The constant $s$ is an invariant depending only on $G$ and $K$. So, it can be computed using a particular case, e.g. when $\Gamma$ is co-compact and $T_n$ have only isolated fixed points. So, Lefschetz's fixed points formula may be used here. We have for example $s=2$ when $G=\PSL(2,\R)$ and $K=\SO(2)$. We can also obtain a speed of convergence in our main theorem in term of the speed of convergence in the equidistribution property of $\widehat T_n$. 
\end{remark}

\section{On the Arnold-Krylov-Guivarc'h theorem} \label{section_rk}

Consider now the case where $G$ is a compact connected semi-simple Lie group, $\Gamma$ is trivial and $K$ a connected compact subgroup of $G$. Define $X:=G/K$. Let $\widehat\lambda$ be the invariant probability measure of $G$ and $\lambda$ its direct image in $X$.

Let $H\subset G$ be a semi-group generated by a finite family 
of elements $g_1,\ldots, g_d$ of $G$.  Denote by $H_n$ the set of words of length $n$ in $H$.
We say that $H$ is {\it equidistributed} on $G$ if for every point $a\in G$, the sequence of probability measures
$$d^{-n}\sum_{g\in H_n} \delta_{ga}$$
converges to $\widehat\lambda$ as $n$ tends to infinity.

The left-multiplication by $g_i$ defines a self-map $\widehat T_{g_i}$ on $G$. Their sum $\widehat T$ can be seen as a correspondence of degree $d$ on $G$. It induces a correspondence $T$ on $X$ of the same degree. So, $H$ is equidistributed if and only if the sequence $\widehat T^n$ is equidistributed. We deduce from our main result the following theorem.

\begin{theorem} \label{th_bis}
Let $G,K,X,\lambda,H$ and $H_n$ be as above. Assume that $H$ is equidistributed on $G$. Then the isolated fixed points in $X$ of the elements of $H_n$  are equidistributed with respect to $\lambda$ when $n$ tends to infinity.
\end{theorem}

Assume that $d=2$ and that the first Betti number of $X$ vanishes.
A result by Guivarc'h \cite{Guivarch} says that if the group generated by $H$ is dense in $G$ then $H$ is equidistributed, see also Arnold-Krylov \cite{ArnoldKrylov}. So, Theorem \ref{th_bis} can be applied in this case.

A similar result holds for groups. Let $H\subset G$ be a group generated by a finite family $\{g_1,\ldots, g_{2d}\}$ where $g_i=g_{2d-i}^{-1}$.  Let $H_n$ denote the family of reduced words of length $n$ in $H$.
We say that $H$ is {\it equidistributed} if the sequence of probability measures
$$\mu_n:={1\over |H_n|}\sum_{g\in H_n} \delta_{ga}$$
converges to $\widehat\lambda$ for every $a\in G$. There are also correspondences $\widehat T_n$ and $T_n$ of degree $|H_n|$ such that $(\widehat T_n)_*(\delta_a)=|H_n|\mu_n$. So, Theorem \ref{th_bis} holds for equidistributed groups $H$.

Another result by Guivarc'h \cite{Guivarch} says that if $d=2$ and if $H$ is dense in $G$ then it is equidistributed. Therefore, our result can be applied under these conditions.

\noindent
T.-C. Dinh, UPMC Univ Paris 06, UMR 7586, Institut de
Math{\'e}matiques de Jussieu, 4 place Jussieu, F-75005 Paris,
France.\\
{\tt  dinh@math.jussieu.fr}, {\tt http://www.math.jussieu.fr/$\sim$dinh}


\small

\begin{thebibliography}{99}

\bibitem{ArnoldKrylov}
Arnold, V. I.; Krylov, A. L.
Uniform distribution of points on a sphere and certain ergodic properties of solutions of linear ordinary differential equations in a complex domain. (Russian) 
{\it Dokl. Akad. Nauk SSSR} {\bf 148} (1963), 9-12. 


\bibitem{BenoistOh}
Benoist, Y.; Oh, H. Equidistribution of rational matrices in their conjugacy classes. {\it Geom. Funct. Anal.} {\bf 17} (2007), no. 1, 1-32.

\bibitem{ClozelOhUllmo}
Clozel, L.; Oh, H.; Ullmo, E.
Hecke operators and equidistribution of Hecke points. 
{\it Invent. Math.} {\bf 144} (2001), no. 2, 327-351. 

\bibitem{ClozelOtal}
Clozel, L.; Otal, J.-P.
Unique ergodicit\'e des correspondances modulaires. {\it Essays on geometry and related topics}, Vol. 1, 2, 205-216, 
Monogr. Enseign. Math., {\bf 38}, {\it Enseignement Math.}, Geneva, 2001.

\bibitem{ClozelUllmo}
Clozel, L.; Ullmo, E. Correspondances modulaires et mesures invariantes. {\it J. Reine Angew. Math.} {\bf 558} (2003), 47-83.

\bibitem{Dinh}
Dinh, T.-C.
Distribution des pr\'eimages et des points p\'eriodiques d'une correspondance polynomiale.  
{\it Bull. Soc. Math. France} {\bf 133} (2005), no. 3, 363-394.

\bibitem{DinhSibony}
Dinh, T.-C.; Sibony, N.
Distribution des valeurs de transformations m\'eromorphes et applications.
{\it Comment. Math. Helv.} {\bf 81} (2006), no. 1, 221-258.
 
\bibitem{Guivarch}
Guivarc'h, Y. G\'en\'eralisation d'un th\'eor\`eme de von Neumann. {\it C. R. Acad. Sci. Paris S\'er. A-B} {\bf 268} (1969), A1020-A1023.
 
\bibitem{Lojasiewicz}
\L ojasiewicz, S.
Sur la g\'eom\'etrie semi- et sous-analytique.
{\it Ann. Inst. Fourier (Grenoble)} {\bf 43} (1993), no. 5, 1575-1595. 

\bibitem{Margulis}
Margulis, G. A.
Discrete subgroups of semisimple Lie groups. 
{\it Ergebnisse der Mathematik und ihrer Grenzgebiete (3)}, {\bf 17}. Springer-Verlag, Berlin, 1991.

\bibitem{Mok}
Mok, N.
Local holomorphic isometric embeddings arising from correspondences in the rank-1 case. {\it Contemporary trends in algebraic geometry and algebraic topology (Tianjin, 2000)}, 155-165, 
Nankai Tracts Math., 5, {\it World Sci. Publ., River Edge, NJ,} 2002.
 
 \bibitem{MokNg}
Mok N.;  Ng S.-C. Germs of measure-preserving holomorphic maps from bounded symmetric domains to their Cartesian products. {\it Preprint} (2010).

\end{thebibliography}
\end{document}